\documentclass[12pt]{article}
 \usepackage{amsmath,amssymb,amsthm}
 \usepackage{color}
 \usepackage{cite}

 \def\draw #1 by #2 (#3){
  \vbox to #2{
    \hrule width #1 height 0pt depth 0pt
    \vfill
    \special{picture #3} 
    }
  }

 \def\scaleddraw #1 by #2 (#3 scaled #4){{
  \dimen0=#1 \dimen1=#2
  \divide\dimen0 by 1000 \multiply\dimen0 by #4
  \divide\dimen1 by 1000 \multiply\dimen1 by #4
  \draw \dimen0 by \dimen1 (#3 scaled #4)}
  }

\newtheorem{theorem}{Theorem}[section]
\newtheorem{example}{Example}
\newtheorem{problem}[example]{Problem}
\newtheorem{defin}[theorem]{Definition}
\newtheorem{lemma}[theorem]{Lemma}

\newtheorem{nt}{Note}

\newenvironment{pf}{\medskip\noindent{Proof:  \hspace*{-.4cm}}
       \enspace}{\hfill \qed \newline \medskip}

 \newcommand{\singlespacing}{\let\CS=\@currsize\renewcommand{\baselinestretch}{1}\tiny\CS}
 \newcommand{\oneandahalfspacing}{\let\CS=\@currsize\renewcommand{\baselinestretch}{1.25}\tiny\CS}
 \newcommand{\doublespacing}{\let\CS=\@currsize\renewcommand{\baselinestretch}{1.35}\tiny\CS}

 \newtheorem{thm}{Theorem}[section]
\newtheorem{cor}{Corollary}[section]
\newtheorem{lem}{Lemma}[section]

\newtheorem{exam}{Example}[section]

\newtheorem{pro}{Proposition}[section]
 \newtheorem{rule-def}[theorem]{Rule}

 \newcommand{\be}{\begin{equation}}
 \newcommand{\ee}{\end{equation}}
 \newcommand{\bea}{\begin{eqnarray}}
 \newcommand{\eea}{\end{eqnarray}}

\def\qed{\hfill \rule{4pt}{7pt}}

\begin{document}

 \begin{center}
 {\Large \bf Nordhaus-Guddum type results for the Steiner Gutman index of graphs}\\

  \vspace{10mm}

 {\bf Zhao Wang}$^{1}$, {\bf Yaping Mao}$^{2,3}$, {\bf Kinkar Chandra Das}$^{4,*},$ {\bf Yilun Shang}$^{5,*}$

 \vspace{9mm}

 \baselineskip=0.20in

 $^1${\it College of Science, China Jiliang University,\\[-1mm]
  Hangzhou 310018, China}\\
 {\rm e-mail:} {\tt wangzhao@mail.bnu.edu.cn}\\[2mm]

 $^2${\it Department of Mathematics, Qinghai Normal University,\\[-1mm]
 Xining, Qinghai 810008, China}\\[2mm]

 $^3${\it Center for Mathematics and Interdisciplinary Sciences of Qinghai Province,\\[-1mm]
 Xining, Qinghai 810008, China}\\
 {\rm e-mail:} {\tt maoyaping@ymail.com}\\[2mm]

 $^4${\it Department of Mathematics, Sungkyunkwan University, \\[-1mm]
 Suwon 16419, Republic of Korea\/} \\
 {\rm e-mail:} {\tt kinkardas2003@googlemail.com}\\[2mm]

 $^4${\it Department of Computer and Information Sciences, Northumbria University, \\[-1mm]
 Newcastle NE1 8ST, UK\/} \\
 {\rm e-mail:} {\tt yilun.shang@northumbria.ac.uk}

 \vspace{4mm}

 (Received September 16, 2020)
 \end{center}

 \baselineskip=0.23in

\begin{abstract}
Building upon the notion of Gutman index $\operatorname{SGut}(G)$,
Mao and Das recently introduced the Steiner Gutman index by
incorporating Steiner distance for a connected graph $G$. The
\emph{Steiner Gutman $k$-index} $\operatorname{SGut}_k(G)$ of $G$ is
defined by $\operatorname{SGut}_k(G)$ $=\sum_{S\subseteq V(G), \
|S|=k}\left(\prod_{v\in S}deg_G(v)\right) d_G(S)$, in which $d_G(S)$
is the Steiner distance of $S$ and $deg_G(v)$ is the degree of $v$
in $G$. In this paper, we derive new sharp upper and lower bounds on
$\operatorname{SGut}_k$, and then investigate the
Nordhaus-Gaddum-type results for the parameter
$\operatorname{SGut}_k$. We obtain sharp upper and lower bounds of
$\operatorname{SGut}_k(G)+\operatorname{SGut}_k(\overline{G})$ and
$\operatorname{SGut}_k(G)\cdot \operatorname{SGut}_k(\overline{G})$
for a connected graph $G$ of order $n$, $m$
edges and maximum degree $\Delta$, minimum degree $\delta$. \\[2mm]
\end{abstract}

 \baselineskip=0.30in

\section{Introduction}

We consider simple, undirected graphs in this paper. The standard
graph-theoretical terminology and notations not defined here follow
\cite{Bondy}. For a graph $G$, let $V(G)$ and $E(G)$ represent its
sets of vertices and edges, respectively. Let $|E(G)|=m$ be the size
of $G$. The complement of $G$ is conventionally denoted by
$\overline{G}$. For a vertex $v\in V(G)$, $deg_G(v)$ is the degree
of $v$. The maximum and minimum degrees are, respectively, denoted
by $\Delta$ and $\delta$. Like degrees, distance is a fundamental
concept of graph theory \cite{BuHa}. For two vertices $u,v\in V(G)$
with connected $G$, the \emph{distance} $d(u,v)=d_G(u,v)$ between
these two vertices is defined as the length of a shortest path
connecting them. An excellent survey paper on this subject can be
found in \cite{Goddard}.

The above classical graph distance has been extended by Chartrand et
al. in 1989 to the Steiner distance, which since then has become an
essential concept of graph theory. Given a graph $G(V,E)$ and a
vertex set $S\subseteq V(G)$ containing no less than two vertices,
\emph{an $S$-Steiner tree} (or \emph{an $S$-tree}, \emph{a Steiner
tree connecting $S$}) is defined as a subgraph $T(V',E')$ of $G$
which is a subtree satisfying $S\subseteq V'$. If $G$ is connected
with order no less than $2$ and $S\subseteq V$ is nonempty, the
\emph{Steiner distance} $d(S)$ among the vertices of $S$ (sometimes
simply put as the distance of $S$) is the minimum size of connected
subgraph whose vertex sets contain the set $S$. Clearly, for a
connected subgraph $H\subseteq G$ with $S\subseteq V(H)$ and
$|E(H)|=d(S)$, $H$ is a tree. When $T$ is subtree of $G$, we have
$d(S)=\min\{|E(T)|\,\,,\,S\subseteq V(T)\}$. For $S=\{u,v\}$,
$d(S)=d(u,v)$ reduces to the classical distance between the two
vertices $u$ and $v$. Another basic observation is that if $|S|=k$,
$d(S)\geq k-1$. For more results regarding varied properties of the
Steiner distance, we refer to the reader to \cite{Ali, Caceresa,
Chartrand, Dankelmann, Goddard, Oellermann}.

In \cite{LMG}, Li et al. generalized the concept of Wiener index
through incorporating the Steiner distance. The \emph{Steiner
$k$-Wiener index} $\operatorname{SW}_k(G)$ of $G$ is defined by
$$
\operatorname{SW}_k(G)=\sum_{\overset{S\subseteq V(G)}{|S|=k}}
d(S)\,.
$$
For $k=2$, it is easy to see the Steiner Wiener index coincides with
the ordinary Wiener index. The interesting range of the Steiner
$k$-Wiener index $\operatorname{SW}_k$ resides in $2 \leq k \leq
n-1$, and the two trivial cases give $\operatorname{SW}_1(G)=0$ and
$\operatorname{SW}_n(G)=n-1$.

Gutman \cite{GutmanSDD} studied the Steiner degree distance, which
is a generalization of ordinary degree distance. Formally, the
\emph{$k$-center Steiner degree distance} $\operatorname{SDD}_k(G)$
of $G$ is given as
$$
\operatorname{SDD}_k(G)=\sum_{\overset{S\subseteq
V(G)}{|S|=k}}\left(\sum_{v\in S}deg_G(v)\right) d_G(S)\,.
$$
\textcolor{red}{The Gutman index of a connected graph $G$ is defined as
$$\operatorname{Gut}(G)=\sum_{u,v\in V(G)}\,deg_G(u)\,deg_G(v)\,d_G(u,v).$$
The Gutman index of graphs attracts attention very recently. For its basic properties and applications, including various lower and upper bounds, see \cite{CL1,DA1,DA2} and the references cited therein.} Recently, Mao and Das \cite{MaoDas} further extended the concept of
Gutman index by incorporating Steiner distance and consider the
weights as a multiplication of degrees. The \emph{Steiner $k$-Gutman
index} $\operatorname{SGut}_k(G)$ of $G$ is defined by
$$
\operatorname{SGut}_k(G)=\sum_{\overset{S\subseteq
V(G)}{|S|=k}}\left(\prod_{v\in S}deg_G(v)\right) d_G(S)\,.
$$
\textcolor{red}{Note that this index is a natural generalization of
classical Gutman index. In particular, for $k=2$, $\operatorname{SGut}_k(G)=Gut(G)$. This is the reason the product of the
degrees comes to the definition of Steiner $k$-Gutman index.} The
weighting of multiplication of degree or expected degree has also
been extensively explored in, for example, the field of random
graphs \cite{ShangH, ShangW} and proves to be very prolific. For
more results on Steiner Wiener index, Steiner degree distance, and
Steiner Gutman index, we refer to the reader to \cite{GutmanSDD,
LMG, MaoDas,MWD1, MWGK, MWGL}.

For a given a graph parameter $f(G)$ and a positive integer $n$, the
well-known \emph{Nordhaus-Gaddum problem} is to determine sharp
bounds for: $(1)$ $f(G)+f(\overline{G})$ and $(2)$ $f(G)\cdot
f(\overline{G})$ over the class of connected graph $G$, with order
$n$, $m$ edges, maximum degree $\Delta$ and minimum degree $\delta$
characterizing the extremal graphs. Many Nordhaus-Gaddum type
relations have attracted considerable attention in graph theory.
Comprehensive results regarding this topic can be found in e.g.
\cite{Aouchiche,Mao,HuaDas,ZhangHu,ShangEI}.

In Section 2, we obtain sharp upper and lower bounds on ${\rm
SGut}_k$ of graph $G$. In Section 3, we obtain sharp upper and lower
bounds of ${\rm SGut}_k(G)+{\rm SGut}_k(\overline{G})$ and ${\rm
SGut}_k(G)\cdot {\rm SGut}_k(\overline{G})$ for a connected graph
$G$ in terms of $n$, $m$, maximum degree $\Delta$ and minimum degree
$\delta$.

\section{Sharp bounds for Steiner Gutman index}

In \cite{MaoDas}, the following results have been obtained:
\begin{lem}{\upshape \cite{MaoDas}}\label{lem2-1}
Let $K_{n}$, $S_n$ and $P_n$ be the complete graph, star graph and
path graph of order $n$, respectively, and let $k$ be an integer
such that $2\leq k\leq n$. Then

$(1)$ ${\rm SGut}_k(K_n)=\binom{n}{k}(n-1)^n(k-1)$;

$(2)$ ${\rm SGut}_k(S_n)=(kn-2k+1)\binom{n-1}{k-1}$;

$(3)$ ${\rm SGut}_k(P_{n})=2^k(k-1)\binom{n}{k+1}$.
\end{lem}

 \noindent
For connected graph $G$ of order $n$ with $m$ edges, the authors in
\cite{MaoDas} derived the following upper and lower bounds on ${\rm
SGut}_k(G)$.
\begin{lem}{\upshape \cite{MaoDas}}\label{lem2-2}
Let $G$ be a connected graph of order $n$ with $m$ edges, and let
$k$ be an integer with $2\leq k\leq n$. Then
$$
(n-1)\left(\frac{2m}{k}\right)^k\binom{n-1}{k-1}^k\geq {\rm
SGut}_k(G)\geq \left\{
\begin{array}{ll}
2m(k-1)\binom{n-1}{k-1} &\mbox {\rm if}~\delta\geq 2\\[2mm]
(k-1)\binom{n}{k}&\mbox {\rm if}~\delta=1.
\end{array}
\right.
$$
\end{lem}

 \noindent
We now give a lower and upper bounds on ${\rm SGut}_k(G)$ in terms
of $n$, $m$, maximum degree $\Delta$ and minimum degree $\delta$:
\begin{pro}\label{pro2-1}
Let $G$ be a connected graph of order $n\geq 3$ with $m$ edges and
maximum degree $\Delta$, minimum degree $\delta$. Also let $k$ be an
integer with $2\leq k\leq n$. Then
\begin{equation*}
2m(n-1){n-1\choose k-1}\frac{\Delta^{k-1}}{k}\geq {\rm
SGut}_k(G)\geq \left\{
\begin{array}{ll}
2m(k-1){n-1\choose
k-1}\frac{\delta^{k-1}}{k}&\mbox {\rm if}~\delta\geq 2\\[4mm]
k{p\choose k}+2^{q}(k-1)\left[{n\choose k}-{p\choose k}\right]&\mbox
{\rm if}~\delta=1,
\end{array}
\right.
\end{equation*}
where $p$ is the number of pendant vertices in $G$, and
$q=\max\{k-p,1\}$. The equality of upper bound holds if and only if
$G$ is a regular graph with $k=n$. The equality of lower bound holds
if and only if $G$ is a regular $(n-k+1)$-connected graph of order
$n$ $(\delta\geq 2)$, or $G\cong P_n$ and $k=n>3$ $(\delta=1)$, or $G\cong P_3$ and $k=2$ $(\delta=1)$.
\end{pro}

\begin{pf} Upper bound: For any $S\subseteq V(G)$ and $|S|=k$, we have $k-1\leq d_G(S)\leq
n-1$, and hence
\begin{equation}
(k-1)\sum_{\overset{S\subseteq V(G)}{|S|=k}}\left(\prod_{v\in
S}deg_G(v)\right)\leq {\rm SGut}_k(G)\leq
(n-1)\sum_{\overset{S\subseteq V(G)}{|S|=k}}\left(\prod_{v\in
S}deg_G(v)\right)\,.\label{1k0}
\end{equation}
Let
$$
M=\sum_{\overset{S\subseteq V(G)}{|S|=k}}\left(\prod_{v\in
S}deg_G(v)\right)=\sum_{\{v_1,v_2,\ldots,v_k\}\subseteq
V(G)}deg_G(v_1)deg_G(v_2)\cdots deg_G(v_k).
$$
and
$$
N=\sum_{\{v_1,v_2,\ldots,v_k\}\subseteq
V(G)}[deg_G(v_1)+deg_G(v_2)+\cdots+deg_G(v_k)].
$$

We first prove the upper bound. Without loss of generality, we can
assume that $deg_G(v_1)\leq deg_G(v_2)\leq \ldots \leq deg_G(v_k)$.
Since
\begin{eqnarray}
&&deg_G(v_1)deg_G(v_2)\ldots deg_G(v_k)\leq \Delta^{k-1}deg_G(v_1)\\[3mm]
&\leq
&\frac{\Delta^{k-1}}{k}(deg_G(v_1)+deg_G(v_2)+\cdots+deg_G(v_k)),\label{1k1}
\end{eqnarray}
it follows that
\begin{eqnarray*}
M&=&\sum_{\{v_1,v_2,\ldots,v_k\}\subseteq V(G)}deg_G(v_1)deg_G(v_2)\ldots deg_G(v_k)\\[3mm]
&\leq &\frac{\Delta^{k-1}}{k}\sum_{\{v_1,v_2,\ldots,v_k\}\subseteq
V(G)}[deg_G(v_1)+deg_G(v_2)+\cdots+deg_G(v_k)]\\[3mm]
&\leq &\frac{\Delta^{k-1}}{k}N.
\end{eqnarray*}
For each $v\in V(G)$, there are ${n-1\choose k-1}$ $k$-subsets in
$G$ such that each of them contains $v$. The contribution of vertex
$v$ is exactly ${n-1\choose k-1}deg_G(v)$. From the arbitrariness of
$v$, we have
$$
N={n-1\choose k-1}\sum_{v\in V(G)}deg_G(v)=2m{n-1\choose k-1},
$$
and hence
 \begin{equation}
{\rm SGut}_k(G)\leq (n-1)M\leq (n-1)\frac{\Delta^{k-1}}{k}N=
2m(n-1){n-1\choose k-1}\frac{\Delta^{k-1}}{k}.\label{1k2}
 \end{equation}

 \vspace*{3mm}

 Suppose that the left equality holds. Then all the inequalities in the
 above must be equalities. From the equality in (\ref{1k1}), one can
 easily see that $G$ is a regular
 graph. From the equality in (\ref{1k2}), we have $d(S)=n-1$ for any $S\subseteq V(G)$, $|S|=k$.
 Since $G$ is connected, then there exists an $S\subseteq V(G)$ such that $|d_G(S)|=k-1$.
 If $k\leq n-1$, then one can easily see that the upper bound is strict as $|d_G(S)|=k-1\leq n-2$
 for some $S$. Otherwise, $k=n$. Since $G$ is connected, we have $|d_G(S)|=n-1$ for any
 $S\subseteq V(G)$. Hence $G$ is a regular graph with $k=n$.

 \noindent
 Conversely, one can see easily that the left equality holds for regular graph with $k=n$.

 \vspace*{3mm}

\noindent
 Lower bound:
Without loss of generality, we can assume that $deg_G(v_1)\leq
deg_G(v_2)\leq \ldots \leq deg_G(v_k)$. First we assume that $\delta\geq 2$. Then
 \begin{eqnarray}
 &&deg_G(v_1)deg_G(v_2)\cdots deg_G(v_k)\geq \delta^{k-1}deg_G(v_k)\nonumber\\[3mm]
 &&~~~~~~~~~~~~~~~~~~~~~~~~~~~~~~~~~~~~~~\geq \frac{\delta^{k-1}}{k}(deg_G(v_1)+deg_G(v_2)+\cdots+deg_G(v_k)),~~~\label{1k3}
 \end{eqnarray}
since $deg_G(v_1)\leq deg_G(v_2)\leq \cdots \leq deg_G(v_k)$.
Furthermore, we have
\begin{eqnarray}
{\rm SGut}_k(G)&\geq&(k-1)\sum_{\{v_1,v_2,\ldots,v_k\}\subseteq V(G)}deg_G(v_1)deg_G(v_2)\ldots deg_G(v_k)~~\label{1k4}\\[3mm]
&\geq&(k-1)\frac{\delta^{k-1}}{k}\sum_{\{v_1,v_2,\ldots,v_k\}\subseteq V(G)}[deg_G(v_1)+deg_G(v_2)+\cdots+deg_G(v_k)]~~\label{1k5}\\[3mm]
&=&(k-1)\frac{\delta^{k-1}}{k}N\nonumber\\
&=&2m(k-1){n-1\choose k-1}\frac{\delta^{k-1}}{k}.\nonumber
\end{eqnarray}

Next we assume that $\delta=1$. If $deg_G(v_1)=deg_G(v_2)=\cdots=deg_G(v_k)=1$,
then $d_G(S)\geq k$ and $deg_G(v_1)deg_G(v_2)\ldots deg_G(v_k)=1$.
If there exists some $v_i$ such that $deg_G(v_i)\geq 2$, then
$d_G(S)\geq k-1$ and $deg_G(v_1)deg_G(v_2)\ldots deg_G(v_k)\geq
2^{\max\{k-p,1\}}=2^{q}$, where $1\leq i\leq k$. Therefore, we have
\begin{eqnarray}
{\rm SGut}_k(G)&\geq&k\sum_{\overset
{\{v_1,v_2,\ldots,v_k\}\subseteq V(G),}{
deg_G(v_1)=deg_G(v_2)=\cdots=deg_G(v_k)=1}}deg_G(v_1)deg_G(v_2)\ldots
deg_G(v_k)\nonumber\\[3mm]
&&+(k-1)\sum_{\overset{\{v_1,v_2,\ldots,v_k\}\subseteq
V(G),}{some~deg_G(v_i)\geq 2}}deg_G(v_1)deg_G(v_2)\ldots deg_G(v_k)~~~\label{1k6}\\
&\geq&k{p\choose k}+2^{q}(k-1)\left[{n\choose k}-{p\choose
k}\right].~~~\label{1k7}
\end{eqnarray}

 \vspace*{3mm}

 Suppose that the right equality holds. Then all the inequalities in the above must be equalities. Suppose that $\delta\geq 2$. From the equality in (\ref{1k4}),
 $d_G(S)=k-1$ for any $S\subseteq V(G)$ and $|S|=k$, that is, $G[S]$ is connected for any $S\subseteq V(G)$ and $|S|=k$, and hence $G$ is
 $(n-k+1)$-connected. From the equality in (\ref{1k5}), we have $deg_G(v_1)=deg_G(v_2)=\cdots=deg_G(v_k)$ for any $S=\{v_1,v_2,\ldots,v_k\}\subseteq V(G)$, and
 hence $G$ is a regular graph. So $G$ is a regular $(n-k+1)$-connected graph of order $n$. 
 
 Next suppose that $\delta=1$. From the equality in (\ref{1k7}), we obtain $deg_G(v_i)=1$ or $deg_G(v_i)=2$ for any vertex $v_i\in V(G)$. Since $G$ is connected, $G\cong P_n$ and $p=2$.
 If $k\geq 3$, then $q=k-p\geq 1$. In this case $d_G(S)=k-1$ for any $S\subseteq V(G)$ and $|S|=k$. One can easily see that $G\cong P_n$ and $k=n>3$ (otherwise, $d_G(S)>k-1$ for some $S\subseteq V(G)$ as $q=k-p$). Otherwise, $k=p=2$ and hence $q=1$. In this case $G\cong P_3$ and $k=2$.

\vspace*{3mm}

 \noindent
 Conversely, one can see easily that the equality holds on lower bound for regular $(n-k+1)$-connected graph of order $n$ $(\delta\geq 2)$, or $G\cong P_n$ and $k=n>3$ $(\delta=1)$, or $G\cong P_3$ and $k=2$ $(\delta=1)$.
 \end{pf}

\begin{exam} \textcolor{red}{Let $G\cong K_n$ with $k=n$. Then 
$$
{\rm SGut}_k(G)=(n-1)^{n+1}=2m(n-1){n-1\choose k-1}\frac{\Delta^{k-1}}{k}.
$$
Let $G\cong K_n\backslash sK_2$ $(n=2s)$ with $k=3$. Then $G$ is a $n-2$ regular graph of order $n$. Then
$$
{\rm SGut}_k(G)=2(n-2)^3{n\choose 3}=2m(k-1){n-1\choose k-1}\frac{\delta^{k-1}}{k}.
$$
Let $G\cong P_n$ with $k=n>3$. Then
$$
{\rm SGut}_k(G)=2^{n-2}(n-1)=k{p\choose k}+2^{q}(k-1)\left[{n\choose k}-{p\choose k}\right]\mbox{ as }p=2.
$$
Let $G\cong P_n$ with $k=2$. Then
$$
{\rm SGut}_k(G)=6=k{p\choose k}+2^{q}(k-1)\left[{n\choose k}-{p\choose k}\right]\mbox{ as }p=2.
$$}
\end{exam}

\section{Nordhaus-Gaddum type results on ${\rm SGut}_k(G)$}

We now in a position to give the Nordhaus-Gaddum type results on
${\rm SGut}_k(G)$.
\begin{thm}\label{th3-2}
Let $G$ be a connected graph of order $n$ with $m$ edges, maximum
degree $\Delta$, minimum degree $\delta$ and a connected
$\overline{G}$. Also let $k$ be an integer with $2\leq k\leq n$.
Then

\noindent
$(1)$
$$
{\rm SGut}_k(G)+{\rm SGut}_k(\overline{G})\leq
(n-1)^2\binom{n}{k}s_1^{k-1}
$$
and
$$
{\rm SGut}_k(G)\cdot {\rm SGut}_k(\overline{G})\leq
2m(n^2-n-2m)(n-1)^2\binom{n-1}{k-1}^2\frac{\Delta^{k-1}\,(n-\delta-1)^{k-1}}{k^2},
$$
where $s_1=\max\{\Delta,\,n-\delta-1\}$. Moreover, the upper bounds
are sharp.

\noindent $(2)$
\begin{eqnarray*}
&&{\rm SGut}_k(G)+{\rm SGut}_k(\overline{G})\\[3mm]
&\geq&\left\{
\begin{array}{ll}
(n-1)(k-1)\binom{n}{k}\,t_1^{k-1} &\mbox {\rm if}~\delta\geq 2, \ \Delta\leq n-3\\[3mm]
2m(k-1){n-1\choose k-1}\frac{\delta^{k-1}}{k}+k{n\choose k}&\mbox {\rm if}~\delta\geq 2, \ \Delta= n-2\\[3mm]
k{n\choose k}+[n(n-1)-2m](k-1){n-1\choose k-1}\frac{(n-\Delta-1)^{k-1}}{k}&\mbox {\rm if}~\delta=1, \ \Delta\leq n-3\\[3mm]
2k{n\choose k}&\mbox {\rm if}~\delta=1, \ \Delta= n-2,
\end{array}
\right.
\end{eqnarray*}
where $t_1=\min\{\delta,\,n-\Delta-1\}$.

$(3)$
\begin{eqnarray*}
&&{\rm SGut}_k(G)\cdot {\rm SGut}_k(\overline{G})\\[3mm]
&\geq&\left\{
\begin{array}{ll}
2m(n^2-n-2m)(k-1)^2\binom{n-1}{k-1}^2\,\frac{\delta^{k-1}\,(n-\Delta-1)^{k-1}}{k^2} &\mbox {\rm if}~\delta\geq 2, \ \Delta\leq n-3\\[3mm]
2m(k-1){n\choose k}{n-1\choose k-1}\delta^{k-1}&\mbox {\rm if}~\delta\geq 2, \ \Delta= n-2\\[3mm]
[n(n-1)-2m](k-1){n\choose k}{n-1\choose
k-1}(n-\Delta-1)^{k-1}&\mbox {\rm if}~\delta=1, \ \Delta\leq n-3\\[3mm]
k^2{n\choose k}^2&\mbox {\rm if}~\delta=1, \ \Delta= n-2.
\end{array}
\right.
\end{eqnarray*}
\end{thm}

\begin{pf} $(1)$ From Proposition \ref{pro2-1}, we have
$$
{\rm SGut}_k(G)\leq 2m(n-1){n-1\choose k-1}\frac{\Delta^{k-1}}{k}
$$
and
$$
{\rm SGut}_k(\overline{G})\leq [n(n-1)-2m](n-1){n-1\choose
k-1}\frac{(n-\delta-1)^{k-1}}{k},
$$
and hence
$$
{\rm SGut}_k(G)+{\rm SGut}_k(\overline{G})\leq
(n-1)^2\binom{n}{k}\,s_1^{k-1}
$$
and
$$
{\rm SGut}_k(G)\cdot {\rm SGut}_k(\overline{G})\leq
2m(n^2-n-2m)(n-1)^2\binom{n-1}{k-1}^2\frac{\Delta^{k-1}(n-\delta-1)^{k-1}}{k^2}.
$$

 \noindent
$(2)$ From Proposition \ref{pro2-1}, if $\delta\geq 2$ and
$\Delta\leq n-3$, then
\begin{eqnarray*}
&&{\rm SGut}_k(G)+{\rm SGut}_k(\overline{G})\\[3mm]
&\geq&2m(k-1){n-1\choose
k-1}\frac{\delta^{k-1}}{k}+[n(n-1)-2m](k-1){n-1\choose
k-1}\frac{(n-\Delta-1)^{k-1}}{k}\\[3mm]
&\geq&(n-1)(k-1)\binom{n}{k}\,t^{k-1}_1.
\end{eqnarray*}

If $\delta(G)\geq 2$ and $\Delta= n-2$, then
\begin{eqnarray*}
&&{\rm SGut}_k(G)+{\rm SGut}_k(\overline{G})\\[3mm]
&\geq&2m(k-1){n-1\choose
k-1}\frac{\delta^{k-1}}{k}+k{p'\choose k}+2^{q'}(k-1)\left[{n\choose k}-{p'\choose k}\right]\\[3mm]
&\geq&2m(k-1){n-1\choose
k-1}\frac{\delta^{k-1}}{k}+k{p'\choose k}+2(k-1)\left[{n\choose k}-{p'\choose k}\right]\\[3mm]
&\geq&2m(k-1){n-1\choose
k-1}\frac{\delta^{k-1}}{k}+k{p'\choose k}+k\left[{n\choose k}-{p'\choose k}\right]\\[3mm]
&=&2m(k-1){n-1\choose k-1}\frac{\delta^{k-1}}{k}+k{n\choose k},
\end{eqnarray*}
where $p'$ is the number of pendant vertices in $G$, and
$q'=\max\{k-p',1\}$.

If $\delta=1$ and $\Delta\leq n-3$, then
\begin{eqnarray*}
&&{\rm SGut}_k(G)+{\rm SGut}_k(\overline{G})\\[3mm]
&\geq&k{p\choose k}+2^{q}(k-1)\left[{n\choose k}-{p\choose
k}\right]+[n(n-1)-2m](k-1){n-1\choose
k-1}\frac{(n-\Delta-1)^{k-1}}{k}\\[3mm]
&\geq&k{n\choose k}+[n(n-1)-2m](k-1){n-1\choose
k-1}\frac{(n-\Delta-1)^{k-1}}{k},
\end{eqnarray*}
where $p$ is the number of pendant vertices in $\overline{G}$, and
$q=\max\{k-p,1\}$.

If $\delta=1$ and $\Delta= n-2$, then
\begin{eqnarray*}
&&{\rm SGut}_k(G)+{\rm SGut}_k(\overline{G})\\[3mm]
&\geq&k{p\choose k}+2^{q}(k-1)\left[{n\choose k}-{p\choose
k}\right]+k{p'\choose k}+2^{q'}(k-1)\left[{n\choose k}-{p'\choose k}\right]\\[3mm]
&\geq&k{n\choose k}+k{n\choose k}\geq 2k{n\choose k},
\end{eqnarray*}
where $p,p'$ are the number of pendant vertices in $G,\overline{G}$,
respectively, and $q=\max\{k-p,1\}$, $q'=\max\{k-p',1\}$.

From the above argument, we have
\begin{eqnarray*}
&&{\rm SGut}_k(G)+{\rm SGut}_k(\overline{G})\\[3mm]
&\geq&\left\{
\begin{array}{ll}
(n-1)(k-1)\binom{n}{k}\,t^{k-1}_1 &\mbox {\rm if}~\delta\geq 2, \ \Delta\leq n-3\\[3mm]
2m(k-1){n-1\choose k-1}\frac{\delta^{k-1}}{k}+k{n\choose
k}&\mbox {\rm if}~\delta\geq 2, \ \Delta= n-2\\[3mm]
k{n\choose k}+[n(n-1)-2m](k-1){n-1\choose
k-1}\,\frac{(n-\Delta-1)^{k-1}}{k}&\mbox {\rm if}~\delta=1, \ \Delta\leq n-3\\[3mm]
2k{n\choose k}&\mbox {\rm if}~\delta=1, \ \Delta= n-2.
\end{array}
\right.
\end{eqnarray*}

\vspace*{3mm}

\noindent For $(3)$, from Proposition \ref{pro2-1}, if $\delta\geq
2$ and $\Delta\leq n-3$, then $${\rm SGut}_k(G)\cdot{\rm
SGut}_k(\overline{G})\geq
2m(n^2-n-2m)(k-1)^2\binom{n-1}{k-1}^2\,\frac{\delta^{k-1}\,(n-\Delta-1)^{k-1}}{k^2}.$$

If $\delta\geq 2$ and $\Delta= n-2$, then
\begin{eqnarray*}
&&{\rm SGut}_k(G)\cdot{\rm SGut}_k(\overline{G})\\[3mm]
&\geq&\left[2m(k-1){n-1\choose
k-1}\frac{\delta^{k-1}}{k}\right]\left[k{p'\choose
k}+2^{q'}(k-1)\left[{n\choose k}-{p'\choose
k}\right]\right]\\[3mm]
&\geq&2m(k-1){n\choose k}{n-1\choose k-1}\delta^{k-1},
\end{eqnarray*}
where $p'$ is the number of pendant vertices in $\overline{G}$, and
$q'=\max\{k-p',1\}$.

If $\delta=1$ and $\Delta\leq n-3$, then
\begin{eqnarray*}
&&{\rm SGut}_k(G)\cdot{\rm SGut}_k(\overline{G})\\[3mm]
&\geq&\left[[n(n-1)-2m](k-1){n-1\choose
k-1}\frac{(n-\Delta-1)^{k-1}}{k}\right]\left[k{p\choose k}+2^{q}(k-1)\left[{n\choose k}-{p\choose k}\right]\right]\\[3mm]
&\geq&[n(n-1)-2m](k-1){n\choose k}{n-1\choose
k-1}\,(n-\Delta-1)^{k-1},
\end{eqnarray*}
where $p$ is the number of pendant vertices in $G$, and
$q=\max\{k-p,1\}$.

If $\delta(G)=1$ and $\Delta= n-2$, then
\begin{eqnarray*}
&&{\rm SGut}_k(G)\cdot{\rm SGut}_k(\overline{G})\\[3mm]
&\geq&\left[k{p\choose k}+2^{q}(k-1)\left[{n\choose k}-{p\choose
k}\right]\right]\left[k{p'\choose k}+2^{q'}(k-1)\left[{n\choose
k}-{p'\choose k}\right]\right]\\[3mm]
&\geq&k^2{n\choose k}^2,
\end{eqnarray*}
where $p,\,p'$ are the number of pendant vertices in $G$ and
$\overline{G}$, respectively, and $q=\max\{k-p,1\}$,
$q'=\max\{k-p',1\}$.

From the above argument, we have
\begin{eqnarray*}
&&{\rm SGut}_k(G)\cdot {\rm SGut}_k(\overline{G})\\[3mm]
&\geq&\left\{
\begin{array}{ll}
2m(n^2-n-2m)(k-1)^2\binom{n-1}{k-1}^2\,\frac{\delta^{k-1}\,(n-\Delta-1)^{k-1}}{k^2} &\mbox {\rm if}~\delta(G)\geq 2, \ \Delta\leq n-3\\[3mm]
2m(k-1){n\choose k}{n-1\choose k-1}\delta^{k-1}&\mbox {\rm if}~\delta(G)\geq 2, \ \Delta= n-2\\[3mm]
[n(n-1)-2m](k-1){n\choose k}{n-1\choose
k-1}(n-\Delta-1)^{k-1}&\mbox {\rm if}~\delta(G)=1, \ \Delta\leq n-3\\[3mm]
k^2{n\choose k}^2&\mbox {\rm if}~\delta(G)=1, \ \Delta= n-2.
\end{array}
\right.
\end{eqnarray*}

\vspace*{3mm}

To show the sharpness of the upper bound and the lower bound for
$\delta(G)\geq 2,\,\Delta\leq n-3$, we let $G$ and $\overline{G}$ be
two $\frac{n-1}{2}$-regular graphs of order $n$, where $n$ is odd.
If $k=n$, then ${\rm SGut}_k(G)=(n-1)(\frac{n-1}{2})^n$, ${\rm
SGut}_k(\overline{G})=(n-1)(\frac{n-1}{2})^n$,
$s_1=\max\{\Delta,\,n-\delta-1\}=\frac{n-1}{2}$,
$\Delta\,(n-\delta-1)=(\frac{n-1}{2})^2$,
$t_1=\min\{\delta,\,n-\Delta-1\}=\frac{n-1}{2}$, and
$\delta\,(n-\Delta-1)=(\frac{n-1}{2})^2$. Furthermore, we have ${\rm
SGut}_k(G)+{\rm SGut}_k(\overline{G})=2(n-1)(\frac{n-1}{2})^n=
(n-1)^2\binom{n}{k}s_1^{k-1}$, ${\rm SGut}_k(G)\cdot {\rm
SGut}_k(\overline{G})=(n-1)^2(\frac{n-1}{2})^{2n}=
2m(n^2-n-2m)(n-1)^2\binom{n-1}{k-1}^2\,\frac{\Delta^{k-1}\,(n-\delta-1)^{k-1}}{k^2}$,
${\rm SGut}_k(G)+{\rm
SGut}_k(\overline{G})=2(n-1)(\frac{n-1}{2})^n=(n-1)(k-1)\binom{n}{k}t_1^{k-1}$
and ${\rm SGut}_k(G)\cdot {\rm
SGut}_k(\overline{G})=(n-1)^2(\frac{n-1}{2})^{2n}=2m(n^2-n-2m)(k-1)^2\binom{n-1}{k-1}^2\,\frac{\delta^{k-1}\,(n-\Delta-1)^{k-1}}{k^2}$.
\end{pf}

\vskip0.3cm

\noindent The following corollary is immediate from above theorem.
\begin{cor}\label{cor4-1}
Let $G$ be a connected graph of order $n\geq 4$ with maximum degree
$\Delta$ and minimum degree $\delta$. Then

\noindent $(1)$
\begin{eqnarray*}
&&(n-1)^2\binom{n}{k}s_1^{k-1}\geq{\rm SGut}_k(G)+{\rm SGut}_k(\overline{G})\\[3mm]
&\geq&\left\{
\begin{array}{ll}
(n-1)(k-1)\binom{n}{k}\,t_1^{k-1} &\mbox {\rm if}~\delta\geq 2, \ \Delta\leq n-3\\[3mm]
n(k-1){n-1\choose k-1}\frac{\delta^{k}}{k}+k{n\choose
k}&\mbox {\rm if}~\delta\geq 2, \ \Delta= n-2\\[3mm]
k{n\choose k}+n(k-1){n-1\choose
k-1}\frac{(n-\Delta-1)^{k}}{k}&\mbox {\rm if}~\delta=1, \ \Delta\leq n-3\\[3mm]
2k{n\choose k}&\mbox {\rm if}~\delta=1, \ \Delta= n-2,
\end{array}
\right.
\end{eqnarray*}
where
$s_1=\min\{\Delta,\,n-\delta-1\}$,\,$t_1=\min\{\delta,\,n-\Delta-1\}$;

\noindent $(2)$
\begin{eqnarray*}
n^2\binom{n-1}{k-1}^2\,\frac{\Delta^{k-1}\,(n-\delta-1)^{k-1}\,(n-1)^4}{4k^2}\geq{\rm SGut}_k(G)\cdot {\rm SGut}_k(\overline{G})\\[3mm]
\geq\left\{
\begin{array}{ll}
n^2(k-1)^2\binom{n-1}{k-1}^2\,\frac{\delta^k\,(n-\Delta-1)^k}{k^2} &\mbox {\rm if}~\delta\geq 2, \ \Delta\leq n-3\\[3mm]
n(k-1){n\choose k}{n-1\choose k-1}\delta^{k}&\mbox {\rm if}~\delta\geq 2, \ \Delta= n-2\\[3mm]
n(k-1){n\choose k}{n-1\choose
k-1}\,(n-\Delta-1)^{k}&\mbox {\rm if}~\delta=1, \ \Delta\leq n-3\\[3mm]
k^2{n\choose k}^2&\mbox {\rm if}~\delta=1, \ \Delta= n-2.
\end{array}
\right.
\end{eqnarray*}
\end{cor}

 \vspace*{3mm}

 \noindent
 The following is the famous inequality by P\'olya and Szeg\"o:
 \begin{lemma} {\rm (P\'olya-Szeg\"o inequality)} {\rm \cite{PS}} \label{1kj1} Let $(a_1,\,a_2,\ldots,\,a_r)$ and $(b_1,\,b_2,\ldots,\,b_r)$ be two positive $r$-tuples such that
 there exist positive numbers $M_1$, $m_1$, $M_2$, $m_2$ satisfying:
 $$0<m_1\leq a_i\leq M_1,~0<m_2\leq b_i\leq M_2,~1\leq i\leq r.$$
 Then
 \begin{equation}
 \frac{\sum\limits^r_{i=1}\,a^2_i\,\sum\limits^r_{i=1}\,b^2_i}{\left(\sum\limits^r_{i=1}\,a_i\,b_i\right)^2}\leq
 \frac{1}{4}\,\left(\sqrt{\frac{M_1\,M_2}{m_1\,m_2}}+\sqrt{\frac{m_1\,m_2}{M_1\,M_2}}\right)^2.\label{1kj2}
 \end{equation}
 \end{lemma}

 \vspace*{3mm}

 \noindent
 We now give another lower and upper bounds on ${\rm SGut}_k(G)\cdot {\rm SGut}_k(\overline{G})$ in terms of $n$, $\Delta$ and $\delta$.
 \begin{thm} Let $G$ be a connected graph of order $n$ with maximum degree $\Delta$, minimum degree $\delta$ and a connected $\overline{G}$. Also let $k$ be an integer with $2\leq k\leq n$. Then
 \begin{equation}
 {\rm SGut}_k(G)\cdot {\rm SGut}_k(\overline{G})\geq\left\{\begin{array}{ll}
 (k-1)^2\,\delta^k\,(n-\delta-1)^k\,\binom{n}{k}^2 &\mbox {\rm if}~\Delta+\delta\leq n-1,\\[3mm]
 (k-1)^2\,\Delta^k\,(n-\Delta-1)^k\,\binom{n}{k}^2 &\mbox {\rm if}~\Delta+\delta\geq n-1
 \end{array} \right.\label{gs1}
 \end{equation}
 with equality holding if and only if $G$ is a regular graph with $d_G(S)=d_{\overline{G}}(S)=k-1$ for any $S\subseteq V(G)$, $|S|=k$,
 and
 $${\rm SGut}_k(G)\cdot {\rm
 SGut}_k(\overline{G})\leq\frac{(n-1)^{2k+2}}{2^{2k+2}}\,\binom{n}{k}^2\,\left[\left(\frac{\Delta\,(n-\delta-1)}{\delta\,(n-\Delta-1)}\right)^{k}+\left(\frac{\delta\,(n-\Delta-1)}{\Delta\,(n-\delta-1)}\right)^{k}+2\right],$$
 Moreover, the equality holds if and only if $G$ is a $\left(\frac{n-1}{2}\right)$-regular graph with $k=n$, $n$ is odd.
 \end{thm}

 \begin{pf} Lower bound: By Cauchy-Schwarz inequality with (\ref{1k0}), we have
  \begin{eqnarray}
 {\rm SGut}_k(G)\cdot {\rm SGut}_k(\overline{G})&\geq& (k-1)^2\,\sum_{\overset{S\subseteq V(G)}{|S|=k}}\left(\prod_{v\in
 S}deg_G(v)\right)\,\sum_{\overset{S\subseteq V(\overline{G})}{|S|=k}}\left(\prod_{v\in S}deg_{\overline{G}}(v)\right)\label{da0}\\
  &\geq&(k-1)^2\,\left(\sum_{\overset{S\subseteq V(G)}{|S|=k}}\left(\prod_{v\in S}deg_G(v)\,\prod_{v\in S}deg_{\overline{G}}(v)\right)^{1/2}\right)^2\label{da1}\\
 &\geq&(k-1)^2\,\left(\sum_{\overset{S\subseteq V(G)}{|S|=k}}\left(\prod_{v\in S}deg_G(v)\,(n-1-deg_{G}(v))\right)^{1/2}\right)^2.\nonumber
 \end{eqnarray}

 \noindent
 Since $\delta\leq deg_G(v)\leq \Delta$, one can easily see that
\begin{eqnarray}
deg_G(v)\,(n-1-deg_{G}(v)) &\geq&\left\{
\begin{array}{ll}
\delta\,(n-\delta-1) &\mbox {\rm if}~\Delta+\delta\leq n-1,\\[3mm]
\Delta\,(n-\Delta-1) &\mbox {\rm if}~\Delta+\delta\geq n-1.
\end{array}\right.\label{1ks1}
\end{eqnarray}

 \noindent
 From the above results, we have
 \begin{eqnarray*}
 {\rm SGut}_k(G)\cdot {\rm SGut}_k(\overline{G})&\geq&\left\{\begin{array}{ll}
 (k-1)^2\,\delta^k\,(n-\delta-1)^k\,\binom{n}{k}^2 &\mbox {\rm if}~\Delta+\delta\leq n-1,\\[3mm]
 (k-1)^2\,\Delta^k\,(n-\Delta-1)^k\,\binom{n}{k}^2 &\mbox {\rm if}~\Delta+\delta\geq n-1.
 \end{array} \right.
 \end{eqnarray*}

 \vspace*{3mm}

 \noindent
 The equality holds in (\ref{da0}) if and only if $d_G(S)=d_{\overline{G}}(S)=k-1$ for any $S\subseteq V(G)$ with $|S|=k$. By Cauchy-Schwarz inequality, the equality holds in (\ref{da1}) if and only if
  $$\frac{\prod_{v\in S_1}deg_G(v)}{\prod_{v\in S_1}deg_{\overline{G}}(v)}=\frac{\prod_{v\in S_2}deg_G(v)}{\prod_{v\in S_2}deg_{\overline{G}}(v)}~\mbox{ for any }~S_1,\,S_2\in V(G)~\mbox{ with }|S_1|=|S_2|=k,$$
 that is, if and only if $deg_G(u)=deg_G(v)$ for any $u,\,v\in V(G)$, that is, if and only if $G$ is a regular graph. Hence the equality holds in (\ref{gs1}) if and only if $G$ is a regular graph with $d_G(S)=d_{\overline{G}}(S)=k-1$ for any $S\subseteq V(G)$, $|S|=k$.

 \vspace*{3mm}

 \noindent
 Upper bound: Let $\overline{\Delta}$ and $\overline{\delta}$ be the maximum degree and the minimum degree of graph $\overline{G}$, respectively. Then $\overline{\Delta}=n-\delta-1$ and
 $\overline{\delta}=n-\Delta-1$. By (\ref{1k0}) and (\ref{1kj2}), we have
 \begin{eqnarray}
 &&{\rm SGut}_k(G)\cdot {\rm SGut}_k(\overline{G})\nonumber\\[3mm]
 &\leq& (n-1)^2\,\sum_{\overset{S\subseteq V(G)}{|S|=k}}\left(\prod_{v\in
 S}deg_G(v)\right)\,\sum_{\overset{S\subseteq V(\overline{G})}{|S|=k}}\left(\prod_{v\in S}deg_{\overline{G}}(v)\right)\nonumber\\
  &\leq&(n-1)^2\,\left(\sum_{\overset{S\subseteq V(G)}{|S|=k}}\left(\prod_{v\in S}deg_G(v)\,\prod_{v\in S}deg_{\overline{G}}(v)\right)^{1/2}\right)^2\,\frac{1}{4}\,\left(\left(\frac{\Delta\,\overline{\Delta}}{\delta\,\overline{\delta}}\right)^{k/2}+\left(\frac{\delta\,\overline{\delta}}{\Delta\,\overline{\Delta}}\right)^{k/2}\right)^2\nonumber\\
 &\leq&\frac{(n-1)^2}{4}\,\left(\sum_{\overset{S\subseteq V(G)}{|S|=k}}\left(\prod_{v\in S}deg_G(v)\,(n-1-deg_{G}(v))\right)^{1/2}\right)^2\,\left(\left(\frac{\Delta\,\overline{\Delta}}{\delta\,\overline{\delta}}\right)^{k/2}+\left(\frac{\delta\,\overline{\delta}}{\Delta\,\overline{\Delta}}\right)^{k/2}\right)^2.\nonumber
 \end{eqnarray}
 One can easily see that
  $$deg_G(v)\,(n-1-deg_{G}(v))\leq \frac{(n-1)^2}{4}~\mbox{ for any }v\in V(G).$$
 Using this result in the above with $\overline{\Delta}=n-\delta-1$ and $\overline{\delta}=n-\Delta-1$, we get
 $${\rm SGut}_k(G)\cdot {\rm
 SGut}_k(\overline{G})\leq\frac{(n-1)^{2k+2}}{2^{2k+2}}\,\binom{n}{k}^2\,\left[\left(\frac{\Delta\,(n-\delta-1)}{\delta\,(n-\Delta-1)}\right)^{k}+\left(\frac{\delta\,(n-\Delta-1)}{\Delta\,(n-\delta-1)}\right)^{k}+2\right].$$

 \vspace*{3mm}

 \noindent
 Moreover, the above equality holds if and only if $G$ is a $\left(\frac{n-1}{2}\right)$-regular graph with $k=n$, $n$ is odd (very similar proof of the Proposition \ref{pro2-1}).
 \end{pf}
 
 \begin{exam} \textcolor{red}{Let $G\cong C_n$ with $k=n$. Then $\delta=2$ and hence
 $$
 {\rm SGut}_k(G)\cdot {\rm SGut}_k(\overline{G})=(n-1)^2(n-3)^n\,2^n=(k-1)^2\,\delta^k\,(n-\delta-1)^k\,\binom{n}{k}^2.
 $$
 Let $G$ be a $\left(\frac{n-1}{2}\right)$-regular graph of order $n$ with $k=n$ and odd $n$. Then $\Delta=\delta=\frac{n-1}{2}$ and hence
 \begin{align*}
 {\rm SGut}_k(G)\cdot {\rm SGut}_k(\overline{G})&=\frac{(n-1)^{2n+2}}{2^{2n}}\\[3mm]
 &=\frac{(n-1)^{2k+2}}{2^{2k+2}}\,\binom{n}{k}^2\,\left[\left(\frac{\Delta\,(n-\delta-1)}{\delta\,(n-\Delta-1)}\right)^{k}+\left(\frac{\delta\,(n-\Delta-1)}{\Delta\,(n-\delta-1)}\right)^{k}+2\right].
 \end{align*}
 }
\end{exam}

 \vspace*{3mm}

\noindent
 We now give another lower and upper bounds on ${\rm SGut}_k(G)+ {\rm SGut}_k(\overline{G})$ in terms of $n$, $\Delta$ and $\delta$.
 \begin{thm} Let $G$ be a connected graph of order $n$ with maximum degree $\Delta$, minimum degree $\delta$ and a connected $\overline{G}$. Also let $k$ be an integer with $2\leq k\leq n$. Then
 \begin{equation}{\rm SGut}_k(G)+{\rm SGut}_k(\overline{G})\geq\left\{ \begin{array}{ll}
 2\,(k-1)\,\delta^{k/2}\,(n-\delta-1)^{k/2}\,\binom{n}{k} &\mbox {\rm if}~\Delta+\delta\leq n-1,\\[5mm]
 2\,(k-1)\,\Delta^{k/2}\,(n-\Delta-1)^{k/2}\,\binom{n}{k} &\mbox {\rm if}~\Delta+\delta\geq n-1
 \end{array}\right.\label{wm1}
 \end{equation}
 with equality holding if and only if $G$ is a $\left(\frac{n-1}{2}\right)$-regular graph with odd $n$ and $d_G(S)=d_{\overline{G}}(S)=k-1$ for any $S\subseteq V(G)$, $|S|=k$, and
 \begin{equation}
 {\rm SGut}_k(G)+{\rm SGut}_k(\overline{G})\leq(n-1)\,\Big[\Delta^k+(n-\delta-1)^k\Big]\,\binom{n}{k}\label{wm2}
 \end{equation}
 with equality holding if and only if $G$ is a regular graph with $k=n$.
 \end{thm}

 \begin{pf} For any two real numbers $a,\,b$, we have $(a-b)^2\geq 0$, that is, $a^2+b^2\geq 2\,a\,b$ with equality holding if and only if $a=b$. Therefore we have
 \begin{eqnarray}
 \prod_{v\in S}deg_G(v)+\prod_{v\in S}deg_{\overline{G}}(v)&\geq& 2\,\left(\prod_{v\in S}deg_G(v)\,\prod_{v\in S}deg_{\overline{G}}(v)\right)^{1/2}\nonumber\\[3mm]
                          &=&2\,\left(\prod_{v\in S}deg_G(v)\,deg_{\overline{G}}(v)\right)^{1/2}\nonumber\\[3mm]
                          &=&2\,\left(\prod_{v\in S}deg_G(v)\,(n-deg_G(v)-1)\right)^{1/2}\,.\nonumber
 \end{eqnarray}

 \noindent
 From the above result with (\ref{1ks1}), we get
 \begin{eqnarray*}
 \prod_{v\in S}deg_G(v)+\prod_{v\in S}deg_{\overline{G}}(v) &\geq&\left\{ \begin{array}{ll}
 2\,\delta^{k/2}\,(n-\delta-1)^{k/2} &\mbox {\rm if}~\Delta+\delta\leq n-1,\\[3mm]
 2\,\Delta^{k/2}\,(n-\Delta-1)^{k/2} &\mbox {\rm if}~\Delta+\delta\geq n-1. \end{array}\right.
 \end{eqnarray*}

 \noindent
 Now,
 \begin{eqnarray*}
 {\rm SGut}_k(G)+{\rm SGut}_k(\overline{G})&=&\sum_{\overset{S\subseteq V(G)}{|S|=k}}\,\left[\left(\prod_{v\in S}deg_G(v)\right) d_G(S)+\left(\prod_{v\in S}deg_{\overline{G}}(v)\right) d_{\overline{G}}(S)\right]\nonumber\\[3mm]
                                          &\geq&(k-1)\,\sum_{\overset{S\subseteq V(G)}{|S|=k}}\,\left[\prod_{v\in S}deg_G(v)+\prod_{v\in S}deg_{\overline{G}}(v)\right]\nonumber\\[3mm]
                                          &\geq&\left\{ \begin{array}{ll}
 2\,(k-1)\,\delta^{k/2}\,(n-\delta-1)^{k/2}\,\binom{n}{k} &\mbox {\rm if}~\Delta+\delta\leq n-1,\\[5mm]
 2\,(k-1)\,\Delta^{k/2}\,(n-\Delta-1)^{k/2}\,\binom{n}{k} &\mbox {\rm if}~\Delta+\delta\geq n-1. \end{array}\right.
 \end{eqnarray*}

 \vspace*{3mm}

 \noindent
 From the above, one can easily see that the equality holds in (\ref{wm1}) if and only if $G$ is a $\left(\frac{n-1}{2}\right)$-regular graph with odd $n$ and $d_G(S)=d_{\overline{G}}(S)=k-1$ for any $S\subseteq V(G)$, $|S|=k$.

 \vspace*{3mm}

 \noindent
 Upper bound: By arithmetic-geometric mean inequality, we have
 \begin{eqnarray*}
 {\rm SGut}_k(G)+{\rm SGut}_k(\overline{G})&=&\sum_{\overset{S\subseteq V(G)}{|S|=k}}\,\left[\left(\prod_{v\in S}deg_G(v)\right) d_G(S)+\left(\prod_{v\in S}deg_{\overline{G}}(v)\right) d_{\overline{G}}(S)\right]\nonumber\\[3mm]
                                          &\leq&(n-1)\,\sum_{\overset{S\subseteq V(G)}{|S|=k}}\,\left[\prod_{v\in S}deg_G(v)+\prod_{v\in S}deg_{\overline{G}}(v)\right]\nonumber\\[3mm]
                                          &\leq&(n-1)\,\sum_{\overset{S\subseteq V(G)}{|S|=k}}\,\left[\left(\frac{\sum\limits_{v\in S}\,deg_G(v)}{k}\right)^k+\left(\frac{\sum\limits_{v\in S}\,deg_{\overline{G}}(v)}{k}\right)^k\right]\nonumber\\[3mm]
                                          &=&\frac{(n-1)}{k^k}\,\sum_{\overset{S\subseteq V(G)}{|S|=k}}\,\left[\left(\sum\limits_{v\in S}\,deg_G(v)\right)^k+\left(\sum\limits_{v\in S}\,(n-deg_G(v)-1)\right)^k\right]\nonumber\\[3mm]
                                          &=&\frac{(n-1)}{k^k}\,\sum_{\overset{S\subseteq V(G)}{|S|=k}}\,\left[\left(\sum\limits_{v\in S}\,deg_G(v)\right)^k+\left(k\,(n-1)-\sum\limits_{v\in S}\,deg_G(v)\right)^k\right]\nonumber\\[3mm]
                                         &\leq&\frac{(n-1)}{k^k}\,\sum_{\overset{S\subseteq V(G)}{|S|=k}}\,\left[(k\,\Delta)^k+\left(k\,(n-1)-k\,\delta\right)^k\right]\nonumber
\end{eqnarray*}
\begin{eqnarray*}
                                          &=&(n-1)\,\Big[\Delta^k+(n-\delta-1)^k\Big]\,\binom{n}{k}.
 \end{eqnarray*}

 \noindent
 From the above, one can easily see that the equality holds in (\ref{wm2})
 if and only if $G$ is a regular graph with $k=n$ (very similar proof of the Proposition \ref{pro2-1}).
 \end{pf}

 \begin{exam} \textcolor{red}{Let $G$ be a $\left(\frac{n-1}{2}\right)$-regular graph with odd $n$ and $k=n$. Then $\delta=\frac{n-1}{2}$ and hence
 $$
 {\rm SGut}_k(G)+{\rm SGut}_k(\overline{G})=\frac{(n-1)^{n+1}}{2^{n-1}}=2\,(k-1)\,\delta^{k/2}\,(n-\delta-1)^{k/2}\,\binom{n}{k}
 $$
 Let $G\cong C_n$ with $k=n$. Then $\Delta=\delta=2$, $\overline{\Delta}=\overline{\delta}=2$ and hence
 $${\rm SGut}_k(G)+{\rm SGut}_k(\overline{G})=(n-1)\,\Big[2^n+(n-3)^n\Big]=(n-1)\,\Big[\Delta^k+(n-\delta-1)^k\Big]\,\binom{n}{k}.$$
 }
\end{exam}


\noindent {\bf Acknowledgment.} \textcolor{red}{The authors are much grateful to three anonymous referees for their valuable comments on our paper, which have
considerably improved the presentation of this paper.} Supported by the National Science
Foundation of China (Nos. 12061059, 11601254, 11661068, 11551001,
11161037, and 11461054), and UoA Flexible Fund from Northumbria
University (201920A1001).

\end{document}